\documentclass[12pt]{article}
\usepackage{amssymb}
\usepackage{amsfonts}
\usepackage{amsmath}
\usepackage{amscd}
\usepackage{amsthm}
\usepackage{latexsym}
\usepackage[all]{xy}
\usepackage{CJK}
\usepackage{setspace}
\usepackage{indentfirst}
\usepackage{extarrows}
\usepackage[dvips]{graphicx}
\usepackage{mathrsfs}
\usepackage{amsmath}
\usepackage{amssymb}
\usepackage{makeidx}
\usepackage{bm}
\usepackage[tight, hang]{subfigure}
\usepackage[dvips]{graphicx}
\usepackage{array}
\usepackage[all]{xy}
\usepackage{mathabx}
\theoremstyle{dinglistyle}
\renewcommand{\eqref}

\usepackage{fancyhdr}

\usepackage{graphicx}
\usepackage{CJK, CJKnumb, amsmath}
\usepackage{mathenv}
\numberwithin{equation}{section}

\usepackage{fancyhdr}
\usepackage{subfigure}

\usepackage{CJK}
\date{}

\date{}

\begin{document}

\title{\bf Relative contravariantly finite subcategories and relative tilting modules $^\star$}
\author{{\small Wei Han, \ Shen Li, \ Shunhua Zhang}\\
         {\small School of Mathematics, Shandong University,
        Jinan, 250100,P.R.China}}
\date{}
\maketitle
\pagenumbering{arabic}
\setlength{\parskip}{0.1\baselineskip}

\begin{abstract}
Let $A$ be a finite dimensional algebra over an algebraically closed field $k$.
Let $T$  be a tilting $A$-module and $B={\rm End}_A\ T$  be the endomorphism algebra of $T$.
In this paper, we consider the correspondence between the tilting $A$-modules and the tilting $B$-modules,
and we prove that there is a one-one correspondence between the basic $T$-tilting $A$-modules in  $T^{\perp}$
and the basic tilting $B$-modules in  $^{\perp}(D_BT)$. Moreover, we  show that  there is a one-one
correspondence between the $T$-contravariantly finite $T$-resolving subcategories of $T^{\perp}$
and the basic $T$-tilting $A$-modules contained in $T^{\perp}$.   As an application,
we show that there is a one-one correspondence between the basic tilting $A$-modules in  $T^{\perp}$
and the basic tilting $B$-modules in $^{\perp}(D_BT)$ if $A$ is a $1$-Gorenstein algebra
or a $m$-replicated algebra over a finite dimensional hereditary algebra.

\end{abstract}

\footnote{ $^\star$ Supported by the National Natural Science Foundation of China (11671230,
11371165, 11171183).}

\footnote{  Email addresses: \  hanwei527@qq.com(H.Han), \   
fbljs603@163.com(S.Li) \\      shzhang@sdu.edu.cn(S.Zhang).}

\vskip0.2in

\section{Introduction}

\vskip 0.2in

Let $A$ be a finite-dimensional algebra over an algebraically closed field $k$.
We denote by mod-$A$ the category of all finitely  generated right $A$-modules,
and we always assume that subcategories of $A$-modules are closed under isomorphisms and direct summands.

\vskip0.2in

Let $\mathcal{M}$ be a subcategory of  $\mbox{mod-A}$. We denote by $\widehat{\mathcal{M}}$
the subcategory of  $\mbox{mod-A}$ consisting of the $A$-modules $L$ such that there is an exact sequence
$0\rightarrow M_{n}\rightarrow M_{n-1}\cdots\rightarrow M_{1}\rightarrow M_{0}\rightarrow
L\rightarrow 0$ with $M_{i}\in\mathcal{M}$, and we define  ${\rm dim}_{\widehat{\mathcal{M}}}(L)$ is the
minimal $n$ such that there exists an exact sequence $0\rightarrow M_{n}\rightarrow M_{n-1}\cdots\rightarrow
M_{1}\rightarrow M_{0}\rightarrow L\rightarrow 0$ with $M_{i}\in\mathcal{M}$.
Dually, we can define $\widecheck{\mathcal{M}}$ and ${\rm dim}_{\widecheck{\mathcal{M}}}(L)$.

\vskip0.2in

For any $0<i<\infty$, the following definition is taken from \cite{AT},
$\mathcal{M}^{\perp_{i}}=\{X\in \mbox{mod}$-$A ~ | ~\mbox{Ext}_{A}^{i}(M, X)=0, \forall M\in\mathcal{M}\}$,
the right orthogonal category of $\mathcal{M}$ is $\mathcal{M}^{\perp}=
\bigcap\limits_{0<i<\infty}\mathcal{M}^{\perp_{i}}=\{ X\in \mbox{mod}$-$A~ | ~ \mbox{Ext}_{A}^{i}(M, X)=0,
\forall  M\in\mathcal{M}, \forall i>0 ~\}$.
Dually, $^{\perp_{i}}\mathcal{M}=\{X\in \mbox{mod}$-$A ~ | ~ \mbox{Ext}_{A}^{i}(X, M)=0,
\forall M\in\mathcal{M}\}$, the left orthogonal category of $\mathcal{M}$ is $^{\perp}\mathcal{M}=
{\bigcap\limits_{0<i<\infty}}^{\perp_{i}}\mathcal{M}=\{X\in \mbox{mod}$-$A ~|~\mbox{Ext}_{A}^{i}(X, M)=0,
\forall i>0, \forall M\in\mathcal{M}\}$. In particular, when $\mathcal{M}={\rm add}M $,
we just denote them by $M^{\perp_{i}}, M^{\perp}, ^{\perp_{i}}M, {^{\perp}M}$ respectively.

\vskip0.2in

Let $T$ be a tilting $A$-module. Denote by ${\cal Y}_T$ the subcategory of ${\rm mod}$-$A$ whose objects
are the $A$-modules $Y$ in $T^{\perp}$ for which there is an exact  sequence
$$ \cdots\xlongrightarrow {f_m} T_{m}\xrightarrow {f_{m-1}} T_{m-1}\xrightarrow{f_{m-2}}\cdots\xlongrightarrow {f_1} T_{1}
\xrightarrow {f_0} T_{0}\longrightarrow Y\longrightarrow 0$$
with $T_{i}$ in ${\rm add}~ T$ and ${\rm Ker} f_i$ in $T^{\perp}$. According to {\cite {AR},
we know that $T^{\perp}={\cal Y}_T$.

\vskip0.2in

Let $M$ be an $A$-module  in  ${\cal Y}_T$. We denote by ${\rm dim}_{\widehat{\mbox{add}T}}(M)$
( $T$-${\rm pd}~ M$ for short) the $T$-projective dimension of $M$, and if $M\in \mbox{add}~ T$,
then $M$ is said to be T-projective. It follows that the subcategory
$\mathcal{T}^{<\infty}(T)=\widehat{\mbox{add}~T}$  consists of all $A$-modules
with finite $T$-projective dimensions. It is easy to see that $\mathcal{T}^{<\infty}(T)$
is the subcategory of $T^{\perp}$. In particular, if $T=A$,
we have  $\mathcal{T}^{<\infty}(T)= \mathcal{P}^{<\infty}(\mbox{mod-A})$.

\vskip 0.2in

Tilting theory is a central topic in the representation theory of algebras,
which has two aspects. One is the external aspect,
which is usually used to compare mod-$A$ to mod-${\rm End}_{A} T$ for a tilting $A$-module $T$, and the
other is the internal aspect, which is to study the structure properties of tilting modules for a fixed
algebra $A$.

\vskip 0.2in

Let $T$ be a tilting $A$-module and $B={\rm End}_A\ T$. Recall from  \cite{M},
Miyashita proved some correspondence between the orthogonal subcategories
of ${\rm mod}$-$A$ and ${\rm mod}$-$B$. In \cite{AR},
Auslander and Reiten have proved  that there is a one-one correspondence
between the basic tilting $A$-modules and the contravariantly finite resolving subcategories of $\mbox{mod-A}$.
However, the relationship between the tilting $A$-modules and the tilting
$B$-modules is little known. In this  paper, we focus on the investigation on the tilting modules in the orthogonal subcategories
and show that there exists a one-one correspondence between the basic $T$-tilting $A$-modules (defined in next section) in $T^{\perp}$
and basic tilting $B$-modules in $^{\bot}(D_{B}T)$.  As an application,
we will prove that there is a one-one correspondence between the basic tilting $A$-modules in  $T^{\perp}$
and the basic tilting $B$-modules in $^{\perp}(D_BT)$ if $A$ is a $1$-Gorenstein algebra
or a $m$-replicated algebra over a finite dimensional hereditary algebra.

\vskip 0.2in

Now, we state our main results in this paper.

\vskip 0.2in

{\bf Theorem 1}\ {\it Let $T$ be a tilting $A$-module and $B={\rm End}_{A}(T)$. Then ${\rm Hom}_{A}(T, -)$ and
$-\otimes _{B}T$ give a one-one correspondence between isomorphism classes of basic $T$-tilting $A$-modules in $T^{\perp}$
and of basic tilting $B$-modules in the subcategory $^{\perp}(D_{B}T)$ of ${\rm mod}$-$B$.}

\vskip 0.2in

We should mention that it is not sure whether a $T$-tilting $A$-module is tilting, but for 1-Gorenstein algebras
 and $m$-replicated algebras over a finite dimensional hereditary algebras, this is true.

\vskip 0.2in

{\bf Theorem 2} {\it Let $A$ be a 1-Gorenstein algebra or a $m$-replicated algebra over a finite dimensional hereditary algebra,
and let $T$ be a tilting $A$-module and $B={\rm End}_{A}(T)$. Then there is a one-one correspondence
between the isomorphism classes of basic tilting $A$-modules in $T^{\perp}$
and of basic tilting $B$-modules in the subcategory $^{\perp}(D_{B}T)$ of ${\rm mod}$-$B$.}

\vskip 0.2in

In general, a $T$-tilting $A$-module is partial tilting. Moreover,we have

\vskip 0.2in

{\bf Theorem 3} \ {\it  Let $T$ be a tilting $A$-module and $B={\rm End}_{A}(T)$. Then ${\rm Hom}_{A}(T, -)$ and
$-\otimes _{B}T$ give a one-one correspondence between basic partial tilting $A$-modules in $T^{\perp}$
and basic partial tilting $B$-modules in $^{\perp}(D_{B}T)$.}

\vskip0.2in

The conceptions of contravariantly and covariantly finite subcategories of ${\rm mod}~A$
were introduced in \cite{AS1,AS2} by Auslander and Smal$\phi$ when they studied
the problem of which subcategories of ${\rm mod}$-$A$ have almost split sequence.
These conceptions have close relationships with tilting modules and cotilting modules, see
\cite{AR} for details.

\vskip0.2in

Let $\mathcal{C}$ be a full subcategory of ${\rm mod}$-$A$,
$C_{M}\in\mathcal{C}$ and $\varphi :C_M\longrightarrow M$ with
$M\in$ $A$-mod. The morphism $\varphi$ is a right
$\mathcal{C}$-approximation of $M$ if the induced  morphism ${\rm
Hom}_A(C,C_{M})\longrightarrow {\rm Hom}_A(C,M)$ is surjective for
any $C\in\mathcal{C}$. A minimal right $\mathcal{C}$-approximation
of $M$ is a right $\mathcal{C}$-approximation which is also a right
minimal morphism, i.e., its restriction to any nonzero summand is
nonzero. The subcategory $\mathcal{C}$ is called contravariantly
finite if any module $M\in$ $A$-mod admits a (minimal) right
$\mathcal{C}$-approximation. The notions of (minimal) left
$\mathcal{C}$-approximation and of covariantly finite subcategory
are dually defined. It is well known that add $M$ is both a
contravariantly finite subcategory and a covariantly finite
subcategory.

\vskip0.2in

Let $\mathcal{X}$ be a subcategory of $\mbox{mod-A}$.
$\mathcal{X}$ is said to be a resolving (resp. coresolving) subcategory
if it is closed under extensions, the kernels of epimorphisms (resp. the cokernels of monomorphisms) and contains all indecomposable projective (resp. injective) $A$-modules. Auslander and Reiten have proved in \cite{AR} that there is a one-one correspondence
between the basic tilting $A$-modules and the contravariantly finite resolving subcategories of $\mbox{mod-A}$.

\vskip0.2in

In this paper, we generalize these kinds of correspondences, and show some
one-one correspondence between the $T$-contravariantly finite $T$-resolving subcategories of $T^{\perp}$
and $T$-tilting $A$-modules contained in $T^{\perp}$.

\vskip0.2in

{\bf Theorem 4}\  {\it Let $T$ be a tilting $A$-module and $B={\rm End}_{A}T$.
Then there is a one-one correspondence between the $T$-resolving $T$-contravariantly finite subcategories of $T^{\perp}$
and the contravariantly finite resolving subcategories  in $^{\perp}(D_{B}T)$.}

\vskip0.2in

Let $T$ be a tilting $A$-module, and let  $L$ and $L^{'}$ be $T$-tilting modules. We say $L$ and $L^{'}$ are equivalent if ${\rm add} L={\rm add} L^{'}$.

\vskip0.2in

{\bf Theorem 5}\ {\it Let $T$ be a tilting $A$-module. Then for any $L\in T^{\perp}$,
 ${\rm add} L\rightarrow \widecheck{\rm add} L\cap T^{\perp}$ and $\mathcal{U}\rightarrow\mathcal{U}\cap\mathcal{U}^{\perp}$ give a one-one correspondence between equivalence class of
basic $T$-tilting $A$-modules in $T^{\perp}$ and the $T$-contravariantly finite $T$-resolving subcategories $\mathcal{U}$ of $T^{\perp}$ which is contained in $\mathcal{T}^{<\infty}(T)$.}

\vskip0.2in

For some special kinds of algebras, we have the following result, which seems to have independent interests.

\vskip0.2in

{\bf Theorem 6} \ {\it Let $A$ be either a 1-Gorenstein algebra or a $m$-replicated
algebra over a hereditary algebra, and $T$ be a tilting $A$-module. Then  there is a one-one correspondence
between the equivalence classes of basic tilting $A$-modules in $T^{\perp}$ and the $T$-contravariantly finite $T$-resolving
subcategories contained in $\mathcal{T}^{<\infty}(T)$.}

\vskip0.2in

This paper is arranged as the following. In section 2, we fix the
notations and recall some necessary facts needed for our further
research. Section 3 is devoted to the proof of Theorem 1, Theorem 2
and Theorem 3. In section 4, we prove Theorem 4, Theorem 5
and Theorem 6.

\vskip0.2in

\section{Preliminary}

Let $A$ be a finite dimensional algebra over an algebraically closed
field $k$. We denote by ${\rm mod}$-$A$ the category of all finitely
generated right $A$-modules and by ${\rm ind}$-$A$  the full subcategory of
${\rm mod}$-$A$ containing exactly one representative of each isomorphism
class of indecomposable $A$-modules. $D={\rm Hom}_k(-,\ k)$ is the
standard duality between $A$-mod and $A^{op}$-mod, and $\tau_A$ is
the Auslander-Reiten translation of $A$.
We denote by ${\rm gl.dim}\ A$ the global dimension of $A$. The Auslander-Reiten quiver
of $A$ is denoted by $\Gamma_A$.

\vskip0.2in

Given an $A$-module $M$, we denote by ${\rm pd} M$ the projective dimension of
$M$ and by ${\rm add} \ M$ the full subcategory
having as objects the direct sums of indecomposable summands of $M$.
For a subcategory  $\mathcal{M}$ of $\mbox{mod}$-$A$, we denote by
$\mbox{add}\mathcal{M}$ the subcategory  of $\mbox{mod}$-$A$
consisting of all direct summand of finitely indecomposable modules in $\mathcal{M}$.

\vskip0.2in

Given any module $M\in A $-mod, we may decompose $M$ as $M\cong\oplus_{i=1}^{m}M_{i}^{d_{i}}$, where
each $M_{i}$ is indecomposable, $d_{i}>0$ for each $i$, and $M_{i}$
is not isomorphic to $M_{j}$ if $i\neq j$. The module $M$ is called
basic if $d_{i}=1 $ for any $i$. The number of non-isomorphic
indecomposable modules occurring in the direct sum decomposition
above is uniquely determined and it is denoted by $\delta(M)$.

\vskip0.2in

An $A$-module $T$ in $\mbox{mod-A}$ is called a (generalized) tilting module if
the following conditions are satisfied:\\
{(1)} ${\rm pd}  T =n<\infty$;\\
{(2)} ${\rm Ext}_A^{i}(T,T)=0 $ for all $i >0$;\\
{(3)} There is  a long exact sequence
 $$0\longrightarrow
A\longrightarrow T_{0}\longrightarrow T_{1}
\longrightarrow\cdot\cdot\cdot\longrightarrow T_{n}\longrightarrow
0$$ with $T_{i}\in {\rm add}\ T $ for $0\leq i\leq n$.

\vskip 0.2in

An $A$-module $M$  satisfying the conditions $(1)$ and $(2)$ of
the definition above is called a partial tilting module.
Let $M$ be a partial tilting module and $X$ be an $A$-module
such that $M\oplus X$ is a tilting module and ${\rm
add}M\cap {\rm add}X=0$. Then X will be called a complement to $M$.

\vskip 0.2in

It is well known that in the classical situation $M$ always admits a
complement and $M$ is a tilting module if and only if $\delta (M) =
\delta (A)$. However, in general situations complements do not
always exist, as shown in \cite{RS}. Moreover it is an important open
problem whether $\delta (M) =
\delta (A)$ is sufficient for a partial tilting module $M$ to be a tilting module.

 \vskip0.2in

{\bf Definition 2.1} {\it Let $\mathcal{U}\subseteq\mathcal{V}\subseteq  {\rm mod} A$. Then we have:

\vskip0.1in

(1)\ If for any $V\in\mathcal{V}$, there is a right $\mathcal{U}$-approximation of V,
then we call $\mathcal{U}$ is contravariantly finite in $\mathcal{V}$.

\vskip0.1in

(2)\ If for any $V\in\mathcal{V}$, there is a left $\mathcal{U}$-approximation of V,
then we call $\mathcal{U}$ is covariantly finite in $\mathcal{V}$.

\vskip0.1in

(3)\ If $\mathcal{U}$ are both contravariantly finite and covariantly finite in $\mathcal{V}$,
we call $\mathcal{U}$ is functorially finite in $\mathcal{V}$.}

\vskip0.2in

Let $T$ be a tilting $A$-module. If a subcategory $\mathcal{U}$ of $T^{\perp}$
is contravariantly finite in $T^{\perp}$, then $\mathcal{U}$
is said to be $T$-contravariantly finite.  Now, we define $T$-resolving subcategory as following.

\vskip0.2in

{\bf Definition 2.2} {\it Let $\mathcal{U}\subseteq T^{\perp}$. Then $\mathcal{U}$
is said to be a $T$-resolving subcategory, if it satisfies the following conditions:

 \vskip0.1in

 (1)\ $\mathcal{U}$ is closed under extensions;

\vskip0.1in

 (2)\  $\mathcal{U}$ is closed under the kernels of epimorphisms in $T^{\perp}$, i.e.,
 for a short exact sequence $0\rightarrow U_1\rightarrow U_2\rightarrow U_3\rightarrow 0$
 in $T^{\perp}$, if $U_2, U_3\in\mathcal{U}$, then $U_2\in \mathcal{U}$.

\vskip0.1in

 (3)\  ${\rm add}T\subseteq\mathcal{U} $.}

 \vskip0.1in

We can define $T$-coresolving subcategory similarly.

\vskip0.2in

{\bf Lemma 2.1}\ {\it  Let $\mathcal{U}$, $\mathcal{V}$ and $\mathcal{W}$
be subcategories of ${\rm mod}$-$A$ with $\mathcal{U}\subseteq\mathcal{V}\subseteq\mathcal{W}$. Then we have:

 \vskip0.1in

  (1) \ If $\mathcal{U}$ is contravariantly finite in $\mathcal{V}$ and  $\mathcal{V}$
  is contravariantly finite in $\mathcal{W}$, then $\mathcal{U}$ is contravariantly finite in $\mathcal{W}$.

 \vskip0.1in

  (2)\ If $\mathcal{U}$ is covariantly finite in $\mathcal{V}$ and $\mathcal{V}$
  is covariantly finite in $\mathcal{W}$, then $\mathcal{U}$ is covariantly finite in $\mathcal{W}$.}

 \vskip0.2in

{\bf Lemma 2.2} {\cite [Proposition 3.7]{AR}}\ {\it  Let $\mathcal{U}$ be a resolving subcategory. Then:

\vskip0.1in

 (1) \ The subcategory of all modules that have right $\mathcal{U}$-approximation is closed under extensions.

 \vskip0.1in

 (2)\ $\mathcal{U}$ is contravariantly finite in ${\rm mod} A$ if and only
 if all simple right A modules have right $\mathcal{U}$-approximations.}

 \vskip0.2in

Let $T$ be a tilting module in $\mbox{mod-A}$, we now define a special partial tilting module, called $T$-tilting module.

 \vskip0.2in

{\bf Defination 2.3} \ {\it Let $T$ be a tilting module in $\mbox{mod-A}$ and $L\in T^{\perp}$.
$L$ is said to be a T-tilting module if it satisfies the following conditions:

\vskip0.1in

(1)  \  ${\rm T}$-${\rm pd} \ L = n < \infty$;

\vskip0.1in

(2)\ ${\rm Ext}_{A}^{i}(L, L)=0$ for $ 0<i<\infty$;

\vskip0.1in

(3) \ There exists an exact sequence $0\rightarrow T\rightarrow L_0\rightarrow L_1
\cdots \rightarrow L_m\rightarrow 0$ with $ L_{i}\in {\rm add} \ L$.}

\vskip0.2in

{\bf Example}\  Let $A=kQ$ be the path algebra of $Q$ with $Q:1\rightarrow 2\rightarrow3$.
The AR quiver $\Gamma_A$ of $A$  is:
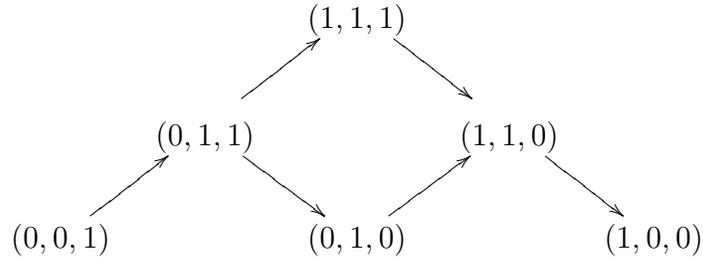
\renewcommand\arraycolsep{0.1cm}
\newcommand{\one}{{\begin{array}{c}\vspace{-5pt}(1, 1, 1)\end{array}}}
\newcommand{\two}{{\begin{array}{c}\vspace{-5pt}(0, 1, 1)\end{array}}}
\newcommand{\three}{{\begin{array}{c}\vspace{-5pt}(1, 1, 0)\end{array}}}
\newcommand{\four}{{\begin{array}{c}\vspace{-5pt}(1, 0, 0)\\\vspace{-5pt}(0, 1, 0)
\\\vspace{-5pt}(0, 0, 1)\\\vspace{-5pt}4\end{array}}}
\[\xymatrix@C=0.3cm@R=0.8cm{
&&\one\ar[dr]&&\\
&\two\ar[ur]\ar[dr]&&\three\ar[dr]&\\
(0, 0, 1)\ar[ur]&&(0, 1, 0)\ar[ur]&&(1, 0, 0)}\]

\vskip0.2in

We denote by $P_{1}=(1, 1, 1), P_{2}=(0, 1, 1), T_{1}=(1, 1, 0) $ and $ T_{2}=(0, 1, 0)$.
Then $T=P_{1}\oplus P_{2}\oplus T_{2}$ is a tilting $A$-module. Let $T^{'}=P_{1}\oplus T_{1}\oplus T_{2}$. Then
 we have an exact sequence $0\rightarrow T\rightarrow P_{1}\oplus P_{1}\oplus T_{2}\oplus T_{2}\rightarrow
T_{1}\rightarrow 0, P_{1}\oplus P_{1}\oplus T_{2}\oplus T_{2}\in \mbox{add}T^{'}$, and $T^{'}$ is a T-tilting module.

\vskip0.2in

Obviously, the $T$-tilting module in the above example is a tilting module,
in particular, $T$-tilting modules over hereditary algebras are tilting modules. In general,
all T-tilting modules are partial tilting modules.

\vskip0.2in

{\bf Proposition 2.3} \ {\it Let $T$ be a tilting $A$-module and  $L$ be a $T$-tilting $A$-module,
then $L$ is a partial tilting $A$-module.}

\vskip0.1in

{\bf Proof}\  Note that $\mbox{Ext}_{A}^{i}(L, L)=0$ for $i>0$ since $L$ is a $T$-tilting module.
 We only need to show $\mbox{pd}_{A}L<\infty$.

Since $T$-pd$L=m<\infty$,  there is an exact sequence $0\rightarrow T_{m}\xlongrightarrow{f_{m}} T_{m-1}
\cdots\xlongrightarrow{f_{1}} T_{0}\xlongrightarrow{f_{0}}
L\longrightarrow 0$ with $T_{i}\in \mbox{add}T$.

Let $C_{i}=\mbox{cokernel}(f_{i})$. Then $L=C_{1}$, and we have following exact sequence:
$$0\longrightarrow T_{m}\longrightarrow T_{m-1} \longrightarrow C_{m}\longrightarrow0$$
$$0\longrightarrow C_{i+1}\longrightarrow T_{i-1} \longrightarrow C_{i}\longrightarrow0 (0<i<m).$$

$T$ is a tilting module, hence $\mbox{pd}_{A}T\leq n<\infty$, and so we get $\mbox{pd}_{A}C_{m}\leq n+1$,
$\mbox{pd}_{A}C_{m-1}\leq n+2$, $\ldots$,
$\mbox{pd}_{A}L=\mbox{pd}_{A}C_{1}\leq n+m<\infty$.  Hence $L$ is a partial tilting $A$-module.      $\hfill\Box$

\vskip0.2in

The following lemmas are well known and useful in our research.

\vskip0.2in

{\bf Lemma 2.4} {\cite [Proposition 3.4 (c)]{AR} \ {\it Let $T$ be a tilting module in ${\rm mod}$-$A$. Then for any $M\in T^{\perp}$,
there is an exact sequence $0\rightarrow K\rightarrow T^{'}\rightarrow M\rightarrow 0$ with $T^{'}\in {\rm add} T$ and $K\in T^{\perp}$.}

\vskip0.2in

{\bf Remark.}\ Let $T$  be a tilting $A$-module. For any $M\in T^{\perp}$,
there exists a long exact sequence $\cdots\rightarrow T_{2}\rightarrow T_{1}
\rightarrow T_{0}\rightarrow M\rightarrow 0$ with $T_{i}\in \mbox{add}T$.

\vskip0.2in

{\bf Lemma 2.5} {\cite [Proposition 1.20]{M}} \ {\it Let $T$ be a tilting $A$-module.
Then for any $M, N\in T^{\perp}$, we have:

\vskip0.1in

(1)\ ${\rm Hom}_{B}({{\rm Hom}}_{A}(T, M), {{\rm Hom}}_{A}(T, N))\cong {{\rm Hom}}_{A}(M, N)$.

\vskip0.1in

(2)\ ${\rm Ext}_{B}^{j}({{\rm Hom}}_{A}(T, M), {{\rm Hom}}_{A}(T, N))\cong {\rm Ext}_{A}^{j}(M, N), 0<j<\infty$.}

\vskip0.2in

{\bf Lemma 2.6} {\cite [Theorem 1.16]{M}}\  {\it Let $T$ be a tilting $A$-module and $B={\rm End}_{A}(T)$.
Then ${{\rm Hom}}_{A}(T, -)$ and
$-\otimes _{B}T$ give a one-one correspondence between the indecomposable modules
in the subcategory $T^{\perp}$ of ${\rm mod}$-$A$ and the indecomposable modules in the subcategory $^{\perp}(D_{B}T)$ of ${\rm mod}$-$B$.}

\vskip0.2in

According to \cite{AR}, we know that there is close relationship between tilting modules
and contravariantly finite (or covariantly finite) subcategories.
Two tilting $A$-modules $T$  and $T'$ are said to be  equivalent if ${\rm add} T={\rm add} T^{'}$.
The following lemma is taken from {\cite{AR}}.

\vskip0.2in

{\bf Lemma 2.7}\ {\cite [Theorem 5.5]{AR}}\ {\it Let $T$  be an $A$-module. Then we have:

 \vskip0.1in

 (1) \  $T\rightarrow T^{\perp}$  gives  a one-one correspondence
 between equivalence class of basic tilting $A$-modules and covariantly
 finite coresolving subcategories $\mathcal{Y}$ with  $\widecheck{\mathcal{Y}}={\rm mod}$-$A$,
 and the inverse correspondence is given by $\mathcal{Y}\rightarrow   {^{\perp}\mathcal{Y}}\cap\mathcal{Y}$.

\vskip0.1in

(2)  \ $T\rightarrow \widecheck{{\rm add} T}$ gives a one-one correspondence between equivalence class of basic tilting modules
and contravariantly finite resolving subcategories $\mathcal{X}$ with  $\mathcal{X}\subseteq
\mathcal{P}^{<\infty}({\rm mod} A)$, and the inverse correspondence is given by  $\mathcal{X}\rightarrow \mathcal{X}\cap\mathcal{X}^{\perp}$.

 \vskip0.1in

 (3)  \  $T\rightarrow {^{\perp}T}$ gives a one-one correspondence between equivalence class of basic cotilting modules
 and contravariantly finite resolving subcategories $\mathcal{X}$ with $\widehat{\mathcal{X}}={\rm mod}$-$A$,
 and the inverse correspondence is given by $\mathcal{X}\rightarrow \mathcal{X}\cap\mathcal{X}^{\perp}$.

 \vskip0.1in

 (4)\  $T\rightarrow \widehat{{\rm add} T}$ gives a one-one correspondence between equivalence class of basic cotilting modules
 and covariantly finite coresolving subcategories $\mathcal{Y}$ with $\mathcal{Y}\subseteq
\mathcal{I}^{<\infty}({\rm mod} A)=\{M\in {\rm mod}$-$A\ |\ {\rm id}_{A}M<\infty\}$,
and the inverse correspondence is given by   $\mathcal{Y}\rightarrow \mathcal{Y}\cap\mathcal{Y}^{\perp}$.}

\vskip0.2in

Throughout this paper, we follow the standard terminology and
notation used in the representation theory of algebras, see \cite{ARS, ASS}

\vskip0.2in

\section{$T$-tilting modules}

 \vskip0.2in

Let $T$  be a tilting $A$-module and $B=\mbox{End}_{A}(T)$. According to Lemma 2.6,
${\mbox{Hom}}_{A}(T, -)$ and $-\otimes _{B}T$ induce a one-one correspondence
between the modules in $T^{\perp}$ and the modules in $^{\perp}(D_{B}T)$.
In this section, we shall show that there exists a one-one correspondence
between the basic $T$-tilting $A$-modules in $T^{\perp}$ and the basic tilting $B$-modules in  $^{\perp}(D_{B}T)$. Moreover, ${\mbox{Hom}}_{A}(T, -)$
and $-\otimes _{B}T$ also give a one-one correspondence between the basic partial tilting $A$-modules
in $T^{\perp}$ and the basic partial tilting $B$-modules in  $^{\perp}(D_{B}T)$.

\vskip0.2in

{\bf Lemma 3.1} \ {\it Let $T$  be a tilting $A$-module and $B={\rm End}_{A}(T)$. Let $M\in T^{\perp}$.
Then ${\rm T}$-${\rm pd} M < \infty$ if and only if ${\rm pd}_{B}({{\rm Hom}}_{A}(T, M)) < \infty$.
In fact, we have ${\rm T}$-${\rm pd} M = {\rm pd}_{B}{\rm  Hom}_{A}(T, M)$.}

\vskip0.1in

{\bf Proof} \  Let $M\in T^{\perp}$ and T-pd$M=k<\infty$. Then there is a long exact sequence:
$$0\rightarrow T_{k}\xlongrightarrow{f_{k}} T_{k-1} \cdots\xlongrightarrow{f_{1}} T_{0}
\xlongrightarrow{f_{0}} M\rightarrow 0, $$
with $T_{i}\in \mbox{add}T\subseteq T^{\perp}$ and $C_{i}=\mbox{cokernel}(f_{i})\in T^{\perp}$.
Applying ${\mbox{Hom}}_{A}(T, -)$ yields an exact sequence:
$$0\rightarrow {\mbox{Hom}}_{A}(T, T_{k})\rightarrow {\mbox{Hom}}_{A}(T, T_{k-1})\cdots
\rightarrow {\mbox{Hom}}_{A}(T, T_{0})\rightarrow {\mbox{Hom}}_{A}(T, M)\rightarrow 0.$$
Note  that $T_{i}\in \mbox{add}T, {\mbox{Hom}}_{A}(T, T)=B$, hence ${\mbox{Hom}}_{A}(T, T_{i})$
are projective $B$-modules, and  $\mbox{pd}_{B}({\mbox{Hom}}_{A}(T, M))\leq k<\infty$.

\vskip0.1in

On the other hand,  we assume that $\mbox{pd}_{B}({\mbox{Hom}}_{A}(T, M))=k<\infty$.
Then there is a projective resolution of ${\mbox{Hom}}_{A}(T, M)$ in $\mbox{mod-B}$:
$$(*) \ \ \   0\rightarrow P_{k}\xlongrightarrow{g_{k}} P_{k-1} \cdots\xlongrightarrow{g_{1}} P_{0}
\xlongrightarrow{g_{0}} {\mbox{Hom}}_{A}(T, M)\rightarrow 0, P_{i}\in
\mbox{add} B.$$
Note that ${\mbox{Hom}}_{A}(T, M)\in {^{\perp}(D_{B}T})$, and $K_{i}=\ker g_{i}\in ^{\perp}(D_{B}T)$. Hence we have
$$D \mbox{Tor}_{1}^{B}(K_i, _{B}T)\cong \mbox{Ext}^{1}_{B}(K_i, {^{\perp}(D_{B}T)})=0.$$
Applying $-\otimes _{B}T$ to $(*)$ yields an exact sequence:
$$0\rightarrow P_{k}\otimes {_{B}T}\rightarrow P_{k-1}\otimes {_{B}T}\cdots\rightarrow P_{0}
\otimes{_{B}T}\rightarrow {\mbox{Hom}}_{A}(T, M)\otimes
{_{B}T}\rightarrow 0 .$$
Since $P_{i}\otimes _{B}T\in $\mbox{add}T and  ${\mbox{Hom}}_{A}(T, M)\otimes _{B}T\cong M$,
we have $T$-${\rm pd}M\leq k<\infty$.

\vskip0.1in

Following from the above proof, we obtain  that $T$-pd$M $ is infinite if and only if $\mbox{pd}_{B}({\mbox{Hom}}_{A}(T, M))$ is infinite.
The proof is completed.      $\hfill\Box$

\vskip0.2in

Let $T$ be a tilting $A$-module and $B=\mbox{End}_{A}(T)$. According to Lemma 2.6, ${\mbox{Hom}}_{A}(T, -)$ and $-\otimes _{B}T$ give a one-one correspondence
between $A$-modules in $T^{\perp}$ and $B$-modules in $^{\perp}(D_{B}T)$. Moreover, we shall show that
they also induce a one-one correspondence between  $T$-tilting $A$-modules in $T^{\perp}$ and  tilting $B$-modules in the subcategory $^{\perp}(D_{B}T)$ of $\mbox{mod-}B$.

\vskip0.2in

{\bf Theorem 3.2}\ {\it Let $T$ be a tilting $A$-module and $B={\rm End}_{A}(T)$. Then ${{\rm Hom}}_{A}(T, -)$ and
$-\otimes _{B}T$ give a one-one correspondence between isomorphism classes of basic $T$-tilting $A$-modules in $T^{\perp}$
and of basic tilting $B$-modules in the subcategory $^{\perp}(D_{B}T)$ of ${\rm mod}$-$B$.}

\vskip0.1in

{\bf Proof} \ \  Let $L\in T^{\perp}$ be a basic $T$-tilting $A$-module. According to Lemma 2.6,  ${\mbox{Hom}}_{A}(T, L)$ is a basic $B$-module belonging to $^{\perp}(D_{B}T)$.
By using Lemma 3.1, we know that $$\mbox{pd}_{B}({\mbox{Hom}}_{A}(T, L))=T-{\rm pd} L<\infty.$$
By Lemma 2.5(2), we know that $\mbox{Ext}_{B}^{j}({\mbox{Hom}}_{A}(T, L), {\mbox{Hom}}_{A}(T, L))\cong\mbox{Ext}_{A}^{j}(L, L)=0,$ for all $j>0$.
Since $L$ is a $T$-tilting $A$-module, we have an exact sequence
$$0\rightarrow T\rightarrow L_0\rightarrow L_1 \cdots \rightarrow L_m\rightarrow 0 \ {\rm with}\  L_{i}\in  {\rm add}\ L.$$
Applying the functor ${\mbox{Hom}}_{A}(T, -)$ yields an exact sequence
$$0\rightarrow {\mbox{Hom}}_{A}(T, T)\rightarrow {\mbox{Hom}}_{A}(T, L_0)\rightarrow {\mbox{Hom}}_{A}(T, L_1) \cdots \rightarrow {\mbox{Hom}}_{A}(T, L_m)\rightarrow 0.$$
It follows that ${\mbox{Hom}}_{A}(T, L)$ is a tilting $B$-module in $^{\perp}(D_{B}T)$.

\vskip0.1in

On the other hand,  let $ M \in {^{\perp}(D_{B}T)}$ be  a basic  tilting  $B$-module. Then $M\otimes _{B}T \in T^{\perp}$
and there  exists an $A$-module $L\in T^{\perp}$ such that $M={\rm Hom}_A(T,L)$  by using Lemma 2.6.

Let ${\rm pd}\ M=m<\infty$. Then there exists a minimal projective resolution of $M$
$$
0\rightarrow P_m\rightarrow \cdots\rightarrow P_1\rightarrow P_0\rightarrow M\rightarrow 0
$$
with $P_i={\rm Hom}_A(T, T_i)$ for some $T_i\in {\rm add}\ T$.  Note that ${\rm Tor}_i^B(M, T)=0$ since
$D{\rm Tor}_i^B(M, T)\simeq {\rm Ext}^i_B(M,DT)=0$ for all $i>0$.
Applying the functor $-\otimes_BT$ to the above sequence yields an exact sequence
$$
0\rightarrow T_m\rightarrow \cdots\rightarrow T_1\rightarrow T_0\rightarrow L\rightarrow 0.
$$
That is $T-{\rm pd}\ L\leq m$.

According to Lemma 2.5 (2),  ${\rm Ext}_A^i(L,L)\simeq {\rm Ext}_B^i({\rm Hom}_A(T,L), {\rm Hom}_A(T,L)=0$.

Finally, since we have an exact sequence
$$
0\rightarrow B\rightarrow M_1\rightarrow M_2\cdots\rightarrow M_s\rightarrow  0
$$
with $M_i\in {\rm add}\ M$ and ${\rm Tor}_i^B(M, T)=0$, applying the functor $-\otimes_BT$ yields an exact sequence
$$
0\rightarrow T\rightarrow L_1\rightarrow L_2\cdots\rightarrow L_s\rightarrow   0,
$$
which means that $M\otimes _{B}T=L$ is a $T$-tilting $A$-module.  The proof is completed.   $\hfill\Box$

\vskip0.2in

${\bf Remark}$\  Let $T$ be a tilting $A$-module and $B=\mbox{End}_{A}(T)$. Let $L$ be a basic
$T$-tilting $A$-module in $T^{\perp}$. It follows that the number of indecomposable direct summands of $L$
is equal to the number of simple $A$-modules. According to Proposition 2.3,  we know that $L$ is also a partial tilting $A$-module.
In general, we are not sure whether a $T$-tilting $A$-module in $T^{\perp}$ is a tilting $A$-module. However, we have the following corollary.

\vskip0.2in

{\bf Corollary 3.3} {\it Let $T$ be a tilting $A$-module and $B={\rm End}_{A}(T)$. Then a $T$-tilting $A$-module is also a tilting  $A$-module, if $A$ is one of the following kinds of algebras:

\vskip0.1in

(1) \ $A$ is a 1-Gorenstein algebra.

\vskip0.1in

(2)\ $A$ is a $m$-replicated algebra over a finite dimensional hereditary algebra.}

\vskip0.1in

{\bf Proof} \ \  One can easily see that if every partial tilting $A$-module $M$ with $\delta (M) = \delta (A)$ is tilting, then
every $T$-tilting $A$-module is also a  tilting  $A$-module.

(1)\ If $A$ is a 1-Gorenstein algebra,  then every partial tilting $A$-module $M$ has projective dimension at most $1$.

\vskip0.1in

(2)\  If $A$ is a $m$-replicated algebra over a finite dimensional hereditary algebra, then according to Theorem 3.1 in \cite{Z},
every partial tilting $A$-module $M$ with $\delta (M) = \delta (A)$ is tilting.         $\hfill\Box$

\vskip0.2in

{\bf Lemma 3.4}\ {\it Let $T$ be a tilting $A$-module and $M\in T^{\perp}$. Then ${\rm pd}_{A}(M)<\infty$
if and only if ${\rm T}$-${\rm pd}\ M < \infty$.}

\vskip0.1in

{\bf  Proof} \  Note that $M\in T^{\perp}$. If T-${\rm pd}  M < \infty$, then by Lemma 2.3, we have $\mbox{pd}_{A}(M)<\infty$.

On the other hand, assume that  $\mbox{pd}_{A}(M)=n<\infty$.According to Lemma 2.4, $M\in T^{\perp}$ implies that there exists an exact sequence:
$$0\rightarrow K_{n}\rightarrow T_{n}\xlongrightarrow{f_{k}} T_{n-1} \cdots\xlongrightarrow{f_{1}} T_{0}\xlongrightarrow{f_{0}} M\rightarrow
0, T_{i}\in \mbox{add}T\subseteq T^{\perp}$$
with $K_{i}=\mbox{ker} f_{i}\in T^{\perp}(0\leq i\leq n)$. It follows that
$$\mbox{Ext}_{A}^{1}(K_{n-1}, K_{n})\simeq\mbox{Ext}_{A}^{2}(K_{n-2}, K_{n})\simeq\cdots \simeq\mbox{Ext}_{A}^{n}(K_{0}, K_{n})\simeq\mbox{Ext}_{A}^{n+1}(M, K_{n})=0$$
since $\mbox{pd}_{A}(M)=n$. Hence the short exact sequence $0\rightarrow K_{n}\rightarrow T_{n}\rightarrow K_{n-1}\rightarrow 0$ splits, and $K_{n-1}$ is a direct summand of $T_{n}$, thus  T-pd$M\leq n<\infty$.  The proof is completed.         $\hfill\Box$

\vskip0.2in

By using Corollary 3.3 and Lemma 3.4,we have

\vskip0.2in

{\bf Theorem 3.5} {\it Let $A$ be a 1-Gorenstein algebra or a $m$-replicated algebra over a finite dimensional hereditary algebra,
and let $T$ be a tilting $A$-module and $B={\rm End}_{A}(T)$. Then there is a one-one correspondence
between the isomorphism classes of basic tilting $A$-modules in $T^{\perp}$
and of basic tilting $B$-modules in the subcategory $^{\perp}(D_{B}T)$ of ${\rm mod} B$.}

\vskip0.1in

{\bf  Proof}\  We only need to show that $T$-tilting $A$-module coincide with tilting $A$-module in $T^{\perp}$.
If $L\in T^{\perp}$ is a $T$-tilting $A$-module, then $L$ is a tilting $A$-module by using Corollary 3.3.

Conversely, assume that $L\in T^{\perp}$ is a tilting $A$-module. By Lemma 3.4, T-${\rm pd}L< \infty$ since ${\rm pd}_A L< \infty$, and
${\rm Ext}^i_A(L,L)=0$ for all $i>0$. Note that whether  $A$ is a 1-Gorenstein algebra or a $m$-replicated algebra,
a tilting $A$-module also is a cotilting module, by using Theorem 5.4 in \cite{AR}, we know that there is an exact sequence
$$
0\rightarrow T\longrightarrow L_{0}\xlongrightarrow{f_{0}} L_{1} \xlongrightarrow{f_{1}} L_{2} \cdots\xlongrightarrow{f_{n}} L_n\rightarrow
0, {\rm with}\  L_{i}\in \mbox{add}T\subseteq T^{\perp}.
$$
Thus $L$ is a $T$-tilting $A$-module.             The proof is completed.         $\hfill\Box$

\vskip0.2in

In general,  we have the following correspondence between partial tilting $A$-modules in $T^{\perp}$
and partial tilting $B$-modules in $^{\perp}(D_{B}T)$.

\vskip0.2in

{\bf Theorem 3.6} \ {\it  Let $T$ be a tilting $A$-module and $B={\rm End}_{A}(T)$. Then ${{\rm Hom}}_{A}(T, -)$ and
$-\otimes _{B}T$ give a one-one correspondence between basic partial tilting $A$-modules in $T^{\perp}$
and basic partial tilting $B$-modules in $^{\perp}(D_{B}T)$.}

\vskip0.1in

{\bf Proof}\  According to Lemma 2.6,  ${\mbox{Hom}}_{A}(T, -)$ and $-\otimes _{B}T$
give a one-one correspondence between basic $A$-modules in $T^{\perp}$  and basic $B$-modules in $^{\perp}(D_{B}T)$.
 By Lemma 2.5(2), for every $A$-module $M$ in $T^{\perp}$, we know that $\mbox{Ext}_{A}^{i}(M, M)=0 (0<i<\infty)$ if and only if  $\mbox{Ext}_{B}^{i}({\mbox{Hom}}_{A}(T, M), {\mbox{Hom}}_{A}(T, M))=0 (0<i<\infty)$. By using Lemma 3.1 and Lemma 3.4,
 for any $M\in T^{\perp}$, we know that $\mbox{pd}_{A}(M)<\infty$ if and only if  $\mbox{pd}_{B}({\mbox{Hom}}_{A}(T, M))<\infty$.
Hence, for an $A$-module $A$ in $T^{\perp}$, $M$ is partial tilting if and only if ${\mbox{Hom}}_{A}(T, M)$ is a tilting $B$-module
belonging to $^{\perp}(D_{B}T)$.  This completes the proof.              $\hfill\Box$

\vskip0.2in

\section{ $T$-contravariantly finite subcategories}

\vskip0.2in

Let $T$ be a tilting $A$-module and $B=\mbox{End}_{A}(T)$. In this section, we prove that there is a one-one correspondence
between basic $T$-tilting $A$-modules in $T^{\perp}$ and the $T$-contravariantly finite $T$-resolving subcategories
in $\mathcal{T}^{<\infty}(T)$.

\vskip0.2in

{\bf Lemma 4.1} \ {\it Let $T$ be a tilting $A$-module and $B={\rm End}_{A}(T)$. Then ${{\rm Hom}}_{A}(T, -)$ and $-\otimes _{B}T$
 induce a one-one correspondence between the $T$-contravariantly finite subcategories of $T^{\perp}$
 and contravariantly finite subcategories of $^{\perp}(D_{B}T)$.}

\vskip0.1in

{\bf Proof}\  \ Let $\mathcal{U}$ be a $T$-contravariantly finite subcategory in $T^{\perp}$. Then
${\mbox{Hom}}_{A}(T, \mathcal{U}) =\{ {\mbox{Hom}}_{A}(T, U)\ | \ U\in\mathcal{U}  \}$ is a subcategory of $^{\perp}(D_{B}T)$.

 For any $Y={\mbox{Hom}}_{A}(T, M)\in $$^{\perp}(D_{B}T)$ with $M\in T^{\perp}$,
 since $\mathcal{U}$ is contravariantly finite in $T^{\perp}$, there exists an $A$-module $U_{M}\in\mathcal{U}$ and
 $f_{M}\in {\mbox{Hom}}_{A}(U_{M}, M)$  such that ${\mbox{Hom}}_{A}(U, f_{M}):{\mbox{Hom}}_{A}(U, U_{M})\rightarrow {\mbox{Hom}}_{A}(U, M)$
 is surjective for any $U\in\mathcal{U}$. Hence $V_{Y}={\mbox{Hom}}_{A}(T, U_{M})\in
{\mbox{Hom}}_{A}(T, \mathcal{U})$, $g_{Y}={\mbox{Hom}}_{A}(T, f_{M})\in {\mbox{Hom}}_{B}(V_{Y}, Y)$.
 For any $V={\mbox{Hom}}_{A}(T, U)\in {\mbox{Hom}}_{A}(T, \mathcal{U})$, we have the following commutative diagram:
\[\xymatrix@C=3cm@R=0.8cm{
   {\mbox{Hom}}_{B}({\mbox{Hom}}_{A}(T, U), {\mbox{Hom}}_{A}(T, U_{M}))\ar[d]_{{\mbox{Hom}}_{B}(V, g_{Y})} \ar[r]^{t_{1}}
                &  {\mbox{Hom}}_{A}(U, U_{M}) \ar[d]^{{\mbox{Hom}}_{A}(U, f_{M})}  \\
   \ar@{}[ur]|-{}  {\mbox{Hom}}_{B}({\mbox{Hom}}_{A}(T, U), {\mbox{Hom}}_{A}(T, M))\ar[r]_{t_{2}}
                &  {\mbox{Hom}}_{A}(U, M)          }\]
where $t_{1}$ and $t_{2}$ are isomorphisms. ${\mbox{Hom}}_{A}(U, f_{M})$ is surjective and this implies that
${\mbox{Hom}}_{B}(V, g_{Y})$ is surjective. Thus ${\mbox{Hom}}_{A}(T, \mathcal{U})$ is contravariantly finite in $^{\perp}(D_{B}T)$.

\vskip0.1in

Similarly, one can prove that if $\mathcal{V}$ is a contravariantly finite subcategory in $^{\perp}(D_{B}T)$,
then $\mathcal{V}\otimes_BT$  is a $T$-contravariantly finite subcategory in $T^{\perp}$.
This  completes the proof.     $\hfill\Box$

\vskip0.2in

{\bf Lemma 4.2} \ {\it Let $T$ be a tilting $A$-module and $B={\rm End}_{A}(T)$.
Then ${{\rm Hom}}_{A}(T, -)$ and $-\otimes _{B}T$ give a one-one correspondence
between the $T$-resolving subcategories of $T^{\perp}$ and the resolving subcategories of $^{\perp}(D_{B}T)$.}

\vskip0.1in

{\bf Proof} \ \   Let $\mathcal{U}$ be a $T$-resolving subcategory of $T^{\perp}$.
Then $\mbox{add}\ B=\mbox{add}\ ({\mbox{Hom}}_{A}(T, T))\subseteq {\mbox{Hom}}_{A}(T, \mathcal{U})$,
and  ${\mbox{Hom}}_{A}(T, \mathcal{U})$ contains all indecomposable projective left $B$-modules.

Suppose  $0\rightarrow V_{1}\rightarrow V_{2}\rightarrow V_{3}\rightarrow 0$ is an exact sequence of $B$-modules
with $V_{1}, V_{3}\in {\mbox{Hom}}_{A}(T, \mathcal{U})\subseteq {^{\perp}(D_{B}T)}$,
then  $V_{2}\in {^{\perp}(D_{B}T)}$ since  $^{\perp}(D_{B}T)$ is closed under extension.
According to Lemma 2.6, there exist $U_{1}$, $U_{3}\in\mathcal{U}$, $U_{2}\in T^{\perp}$
such that $V_{i}={\mbox{Hom}}_{A}(T, U_{i}), U_{i}=V_{i}\otimes_{B}T(i=1, 2, 3)$.
Since  $D\mbox{Tor}_{1}(V_{3}, T)\cong \mbox{Ext}_{B}^{1}(V_{3}, D_{B}T)=0$, applying the functor $-\otimes_{B}T$ to the exact sequence $0\rightarrow
V_{1}\rightarrow V_{2}\rightarrow V_{3}\rightarrow 0$, we obtain an exact sequence
 $0\rightarrow V_{1}\otimes_{B}T\rightarrow V_{2}\otimes_{B}T\rightarrow V_{3}\otimes_{B}T\rightarrow 0$, hence we have
 an exact sequence $0\rightarrow U_{1}\rightarrow U_{2}\rightarrow U_{3}\rightarrow 0$ in $T^{\perp}$.
 Since $\mathcal{U}$ is $T$-resolving in $T^{\perp}$, $U_{1}$ and $U_{3}\in\mathcal{U}$,
 we know that $U_{2}\in\mathcal{U}$, and  $V_{2}={\mbox{Hom}}_{A}(T, U_{2})\in {\mbox{Hom}}_{A}(T, \mathcal{U})$,
i.e., ${\mbox{Hom}}_{A}(T, \mathcal{U})$ is closed under extension.

Similarly, one can prove that ${\mbox{Hom}}_{A}(T, \mathcal{U})$ is closed under the kernel of epimorphism.
Hence ${\mbox{Hom}}_{A}(T, \mathcal{U})$ is a resolving subcategory of $^{\perp}(D_{B}T)$.

\vskip0.1in

On the other hand, by using the same method, we can prove that if $\mathcal{V}$ is a resolving subcategory of $^{\perp}(D_{B}T)$,
then $\mathcal{V}\otimes_{B}T$ is a $T$-resolving subcategory in $T^{\perp}$. The proof is completed.       $\hfill\Box$

\vskip0.2in

Summarizing  Lemma 4.1 and Lemma 4.2, we have the following.

\vskip0.2in

{\bf Theorem 4.3}\  {\it Let $T$ be a tilting $A$-module and $B={\rm End}_{A}T$.
Then there is a one-one correspondence between the $T$-contravariantly finite  $T$-resolving subcategories of $T^{\perp}$
and the contravariantly finite resolving subcategories  in $^{\perp}(D_{B}T)$.}

\vskip0.2in

Let $T$ be a tilting $A$-module and $B=\mbox{End}_{A}T$. According to Lemma 2.7, there is a one-one correspondence
between the equivalence classes of basic tilting $B$-modules and contravariantly finite subcategories of ${\rm mod}$-$B$.
 Since $^{\perp}(D_{B}T)$ is a resolving and contravariantly finite subcategory of ${\rm mod}$-$B$, restricting
 the above one-one correspondence on the subcategory $^{\perp}(D_{B}T)$ also yields the following lemma.

\vskip0.2in

{\bf Lemma 4.4} \ {\it Let $T$ be a tilting $A$-module and $B={\rm End}_{A}T$.  Let $M\in$$^{\perp}(D_{B}T)$.
Then $M\rightarrow \widecheck{{\rm add} M}$ and $\mathcal{X}\rightarrow \mathcal{X}\cap\mathcal{X}^{\perp}$
give a one-one correspondence between the equivalence classes of basic tilting $B$-modules in $^{\perp}(D_{B}T)$
and the resolving contravariantly finite  subcategories
belonging to $\mathcal{X}\subseteq \mathcal{P}^{<\infty}({\rm modB})\cap {^{\perp}(D_{B}T)}$.}

\vskip0.2in

{\bf Lemma 4.5}\ {\it Let $T$ be a tilting $A$-module, $B={\rm End}_{A}T$, and
$\mathcal{X}$ be a subcategory of ${^{\perp}(D_{B}T)}$. Then we have following:

\vskip0.1in

 (1)\  $\widecheck{\mathcal{X}}\otimes_{B} T={\widecheck {\mathcal{X}\otimes_{B} T}}\cap T^{\perp}$.

\vskip0.1in

(2)\ $(\mathcal{X}^{\perp}\cap$$^{\perp}(D_{B}T))\otimes_{B}T=(\mathcal{X}\otimes _{B}T)^{\perp}\cap T^{\perp}$.}

\vskip0.1in

{\bf Proof} \  \  (1)\  Let $M\otimes_{B} T\in\widecheck{\mathcal{X}}\otimes_{B}T$ with $M\in\widecheck{\mathcal{X}}$. Then there is an exact sequence:
$$0\rightarrow M\rightarrow X_{0}\xlongrightarrow{f_{0}} X_{1} \cdots\xlongrightarrow{f_{n-1}} X_{n}\rightarrow 0\ {\rm with}\  X_{i}\in\mathcal{X}. $$
Let $K_{i}=\mbox{ker} f_{i}$ and $M=\mbox{ker} f_{0}=K_{0}$. Since $^{\perp}(D_{B}T)$ is resolving in $\mbox{mod-}B$ and $\mathcal{X}\subseteq$$^{\perp}(D_{B}T)$, we have that $K_{i}\in {^{\perp}(D_{B}T)}$ and $M=K_{0}\in $$^{\perp}(D_{B}T)$, hence $M\otimes _{B}T\in T^{\perp}$.
Note that $D\mbox{Tor}_{1}(K_{i}, _{B}T)\cong\mbox{Ext}_{B}^{1}(K_{i}, D_{B}T)=0$, applying the functor $-\otimes_{B} T$ on the above sequence,
we obtain an exact sequence:
$$0\rightarrow M\otimes_{B} T\rightarrow X_{0}\otimes_{B} T\rightarrow X_{1}\otimes_{B} T \cdots\rightarrow X_{n}\otimes_{B} T\rightarrow 0, $$
it follows that $M\otimes_{B} T\in\widecheck{\mathcal{X}\otimes_{B} T}$, hence $M\otimes_{B} T\in\widecheck{\mathcal{X}\otimes_{B} T}\cap T^{\perp}$.

On the other hand,  for any $N\in\widecheck{\mathcal{X}\otimes_{B}T}\cap T^{\perp}$,
there exists an exact sequence
$$0\rightarrow N\rightarrow Y_{0}\otimes{_{B}T}\rightarrow Y_{1}\otimes{_{B}T} \cdots\rightarrow Y_{m}\otimes{_{B}T}\rightarrow 0$$
 with $Y_i\in\mathcal{X}$.  Applying the functor ${\mbox{Hom}}_{A}(T, -)$, by using Lemma 2.6 and $N\in T^{\perp}$, we have an exact sequence
$0\rightarrow {\mbox{Hom}}_{A}(T, N)\rightarrow Y_{0}\rightarrow Y_{1} \cdots\rightarrow Y_{m}\rightarrow 0.$
Hence ${\mbox{Hom}}_{A}(T, N)\in\widecheck{\mathcal{X}}$, and $N= {\mbox{Hom}}_{A}(T, N)\otimes{_{B}T}\in\widecheck{\mathcal{X}}\otimes _{B}T$.

\vskip0.1in

(2)\ \ Let $M\otimes _{B}T\in(\mathcal{X}^{\perp}\cap$$^{\perp}(D_{B}T))\otimes_{B}T$ with $M\in\mathcal{X}^{\perp}\cap$$^{\perp}(D_{B}T)$.
 According to Lemma 2.5 (2),  for any $X\otimes_{B}T\in\mathcal{X}\otimes _{B}T$ with $X\in\mathcal{X}$, we have
$\mbox{Ext}_{A}^{i}(X\otimes_{B}T, M\otimes_{B}T)\cong \mbox{Ext}_{B}^{i}(X, M)=0,$ hence $M\otimes_{B}T\in(\mathcal{X}\otimes _{B}T)^{\perp}$.
Note that $M\in^{\perp}(D_{B}T)$, thus $M\otimes_{B}T\in T^{\perp}$, hence  $M\otimes_{B}T\in(\mathcal{X}\otimes _{B}T)^{\perp}\cap T^{\perp}$.

On the other hand,  let $N\in(\mathcal{X}\otimes _{B}T)^{\perp}\cap T^{\perp}$. Then ${\mbox{Hom}}_{A}(T, N)\in$$^{\perp}(D_{B}T)$.
For every $X\in\mathcal{X}$, by using  Lemma 2.5 (2) again, we have
$$\mbox{Ext}_{B}^{i}(X, \mbox{Hom}_{A}(T, N))\simeq \mbox{Ext}_{A}^{i}(X\otimes_B T, \mbox{Hom}_{A}(T, N)\otimes_B T)
\simeq \mbox{Ext}^i_{A}(X\otimes_B T, N)=0.$$
It follows that ${\mbox{Hom}}_{A}(T, N)\in \mathcal{X}^{\perp}$ and ${\mbox{Hom}}_{A}(T, N)\in \mathcal{X}^{\perp}\cap$$^{\perp}(D_{B}T)$,
Hence $N={\mbox{Hom}}_{A}(T, N)\otimes_{B} T\in(\mathcal{X}^{\perp}\cap {^{\perp}(D_{B}T))}\otimes_{B} T$.
The proof is completed.      $\hfill\Box$

\vskip0.2in

Let $T$ be a tilting $A$-module, and let  $L$ and $L^{'}$ are $T$-tilting modules. We say $L$ and $L^{'}$ are equivalent if ${\rm add} L={\rm add} L^{'}$.

\vskip0.2in

{\bf Theorem 4.6}\ {\it Let $T$ be a tilting $A$-module. Then for any $L\in T^{\perp}$, 
 ${\rm add} L\rightarrow \widecheck{{\rm add} L}\cap T^{\perp}$ and
$\mathcal{U}\rightarrow\mathcal{U}\cap\mathcal{U}^{\perp}$
give a one-one correspondence between equivalence class of
basic $T$-tilting $A$-modules in $T^{\perp}$ and the
$T$-contravariantly finite $T$-resolving subcategories $\mathcal{U}$ of $T^{\perp}$
which is contained in $\mathcal{T}^{<\infty}(T)$.}

\vskip0.1in

{\bf Proof}\  \  \ For any $X\in T^{\perp}$ and $Y\in$$^{\perp}(D_{B}T)$, it is easy to see that
${\mbox{Hom}}_{A}(T, \mbox{add} X)=\mbox{add} {\mbox{Hom}}_{A}(T, X)$,
$(\mbox{add}Y)\otimes_{B} T=\mbox{\mbox{add}} (Y\otimes_{B} T)$. According to Lemma 3.1 and Lemma 2.6,
we have  T-pd$X=$$\mbox{pd}_{B}({\mbox{Hom}}_{A}(T, X))$ and
${\mbox{Hom}}_{A}(T, \mathcal{T}^{<\infty}(T))=\mathcal{P}^{<\infty}(\mbox{modB})\cap{^{\perp}(D_{B}T)}$.

\vskip0.1in

(1) \ Let $L$ be a basic $T$-tilting $A$-module in $T^{\perp}$. Then by Theorem 3.2,
$Y={\mbox{Hom}}_{A}(T, L)\in {^{\perp}(D_{B}T)}$ is a basic tilting $B$-module.
By Lemma 4.4, we know that $\widecheck{\mbox{add} Y}$ is contravariantly finite in $\mbox{mod-B}$,
which is  contained in $\mathcal{P}^{<\infty}(\mbox{modB})\cap{^{\perp}(D_{B}T)}$.
Applying the functor $-\otimes_{B} T$ on $\widecheck{\mbox{add}Y}$, and using Lemma 4.4 and Lemma 4.5,
we know that  $\widecheck{\mbox{add}Y}\otimes _{B} T=\widecheck{\mbox{add}L}\cap T^{\perp}$ is a $T$-contravariantly finite $T$-resolving subcategory
 of $T^{\perp}$ which is contained in $\mathcal{T}^{<\infty}(T)$.

\vskip0.1in

(2)  \   Let $\mathcal{U}$ be a $T$-contravariantly finite $T$-resolving subcategory of $T^{\perp}$
and $\mathcal{U}\subseteq \mathcal{T}^{<\infty}(T)$.
Then ${\mbox{Hom}}_{A}(T, \mathcal{U})\subseteq\mathcal{P}^{<\infty}(\mbox{modB})\cap{^{\perp}(D_{B}T)}$.
According to Lemma 4.1 and Lemma 4.2,  we know that ${\mbox{Hom}}_{A}(T, \mathcal{U})$
is a contravariantly finite resolving  subcategory, which is contained in $^{\perp}(D_{B}T)$.
By using Lemma 4.4,  there exists a tilting $B$-module $M\in^{\perp}(D_{B}T)$  such that
${\mbox{Hom}}_{A}(T, \mathcal{U})\cap ({\mbox{Hom}}_{A}(T, \mathcal{U}))^{\perp}=\mbox{add} M$,
hence  ${\mbox{Hom}}_{A}(T, \mathcal{U})\cap(({\mbox{Hom}}_{A}(T, \mathcal{U}))^{\perp}\cap$ $^{\perp}(D_{B} T))=\mbox{add}M$.
Applying the functor $-\otimes _{B}T$ and using Lemma 4.5,
we have $\mathcal{U}\cap (\mathcal{U}^{\perp}\cap T^{\perp})=\mbox{add}(M\otimes_{B}T)$, thus
$\mathcal{U}\cap \mathcal{U}^{\perp} =\mbox{add}(M\otimes_{B}T)$,
since $\mathcal{U}$ is a subcategory of $T^{\perp}$. According to Theorem 3.2,
$M\otimes_{B} T$ is a $T$-tilting $A$-module in $T^{\perp}$.

\vskip0.1in

(3) \   Now we show that $\mbox{add}L\rightarrow \widecheck{\mbox{add}L}\cap T^{\perp}$
and $\mathcal{U}\rightarrow\mathcal{U}\cap\mathcal{U}^{\perp}$ are a pair of converse functors.

Let $L$ be a basic $T$-tilting $A$-module. Then $Y={\mbox{Hom}}_{A}(T, L)\in {^{\perp}(D_{B}T)}$
is a basic tilting $B$-module. According to Lemma 4.4, we have that
$\widecheck{\mbox{add}Y}\cap(\widecheck{\mbox{add}Y})^{\perp}=\mbox{add}Y$ and
$\mbox{add}Y=\widecheck{\mbox{add}Y}\cap(\widecheck{\mbox{add}Y})^{\perp}\cap ^{\perp}(D_{B}T)$.
 Applying the the functor $-\otimes _{B}T$ and using  Lemma 4.5,
 we have
 $$\mbox{add}L=(\widecheck{\mbox{add}L}\cap T^{\perp})\cap (\widecheck{\mbox{add}L}\cap T^{\perp})^{\perp}\cap T^{\perp}
 =(\widecheck{\mbox{add}L}\cap T^{\perp})\cap (\widecheck{\mbox{add}L}\cap T^{\perp})^{\perp}.$$

 On the other hand, let $\mathcal{U}$ be a $T$-contravariantly finite $T$-resolving subcategory
 and $\mathcal{U}\subseteq \mathcal{T}^{<\infty}(T)$.
 By Lemma 4.2,  $\mathcal{V}={\mbox{Hom}}_{A}(T, \mathcal{U})$ is a contravariantly finite resolving subcategory,
 which is contained in $\mathcal{P}^{<\infty}(\mbox{modB})\cap{^{\perp}(D_{B}T)}$.
 By using Lemma 4.4, we have $\widecheck{\mathcal{V}\cap\mathcal{V}^{\perp}}=\mathcal{V}$.
 Applying the functor $-\otimes _{B}T$ on the both side and using Lemma 4.5,
 we have $\mathcal{U}=\widecheck{\mathcal{U}\cap\mathcal{U}^{\perp}}\cap T^{\perp}$.
The proof is finished.            $\hfill\Box$

\vskip0.2in

Combining Theorem 3.5 and Theorem 4.6, we obtain the following result, which seems to have independent interests.

\vskip0.2in

{\bf Theorem 4.7} \ {\it Let $A$ be either a 1-Gorenstein algebra or a $m$-replicated
algebra over a hereditary algebra, and $T$ be a tilting $A$-module. Then  there is a one-one correspondence
between the equivalence classes of basic tilting $A$-modules in $T^{\perp}$ and the $T$-contravariantly finite $T$-resolving
subcategories contained in $\mathcal{T}^{<\infty}(T)$.}

\end{document}